\documentclass[12pt]{article}
\usepackage[utf8]{inputenc}
\usepackage[english]{babel}
\usepackage{amsmath,amsfonts,amssymb,amsthm}
\usepackage{thmtools}
\usepackage{bm}
\usepackage{geometry}
\usepackage{mathtools}
\usepackage{nomencl}
\usepackage{tikz,pgfplots}
\usepackage{caption}
\usepackage{rotating}
\usepackage[labelfont={rm},position=below]{subfig}
\usepackage{setspace}
\usepackage{subfig}
\usepackage[small]{titlesec}
\usepackage{hyperref}

\allowdisplaybreaks 

\pgfplotsset{compat=1.16}
\usetikzlibrary{
  angles, arrows.meta, calc, decorations.markings, decorations.pathmorphing,
  fadings, patterns, quotes, intersections
}

\DeclareCaptionFormat{notext}{textformat=empty}

\hypersetup{
  colorlinks=false,
  pdfborderstyle={/S/U/W 0.5},
  linkbordercolor={1 0 0},
  pdftitle={},
  pdfauthor={Keegan Doig Anderson}
}

\newcommand{\beq}{\begin{equation}}
\newcommand{\eeq}{\end{equation}}

\newcommand{\R}{\mathbb{R}}

\newcommand{\defeq}{\coloneqq} 

\title{Dynamics of a Rotated Orthogonal Gravitational Wedge Billiard}
\author{K. D. Anderson}
\date{}

\begin{document}
\maketitle

\begin{abstract}
  We investigate a rotated, orthogonal gravitational wedge billiard---a special
  case of the asymmetric gravitational wedge billiard---in which the dynamics are
  integrable.
  We derive equations and conditions under which periodic orbits may be
  constructed for this model, and show that any other trajectory will be dense
  in the configuration space.
\end{abstract}


\section{Introduction}
Based on previous work on the asymmetric wedge billiard
\cite{anderson2019thesis,anderson2019computational}, we now investigate
a special case of the asymmetric wedge billiard which leads to integrable
dynamics.

\begin{figure}
  \begin{center}
    \includegraphics{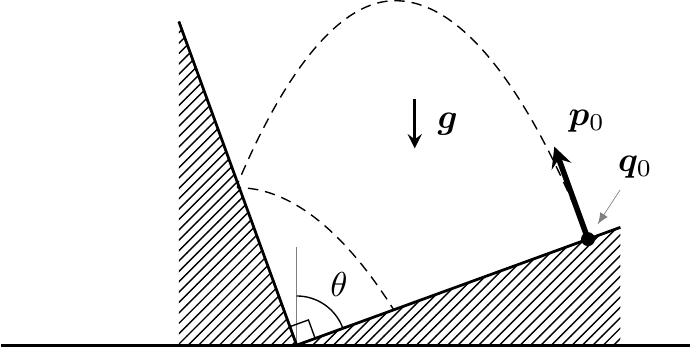}
    \caption{The rotated orthogonal wedge billiard.}
    \label{fig:awb-orthogonal-wedge}
  \end{center}
\end{figure}

For this special case, we set $\theta_1 = \theta$ and $\theta_2 = \pi/2 -
\theta$ in the asymmetric wedge billiard
\cite{anderson2019thesis,anderson2019computational}.
This corresponds to an orthogonal wedge which is rotated by an angle $\theta$
from the vertical, see Figure \ref{fig:awb-orthogonal-wedge}.

\begin{figure}[ht]
    \begin{center}
    \subfloat[\label{fig:owb-fixed-point-trajectory}] {
      \includegraphics[width=.4\textwidth]{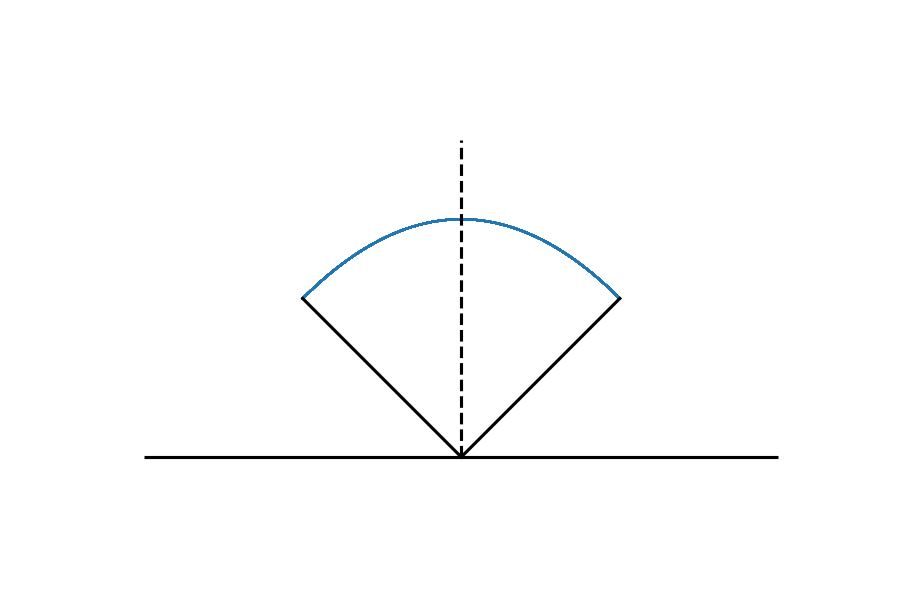}
    }
    \subfloat[\label{fig:owb-period-two-trajectory}] {
      \includegraphics[width=.4\textwidth]{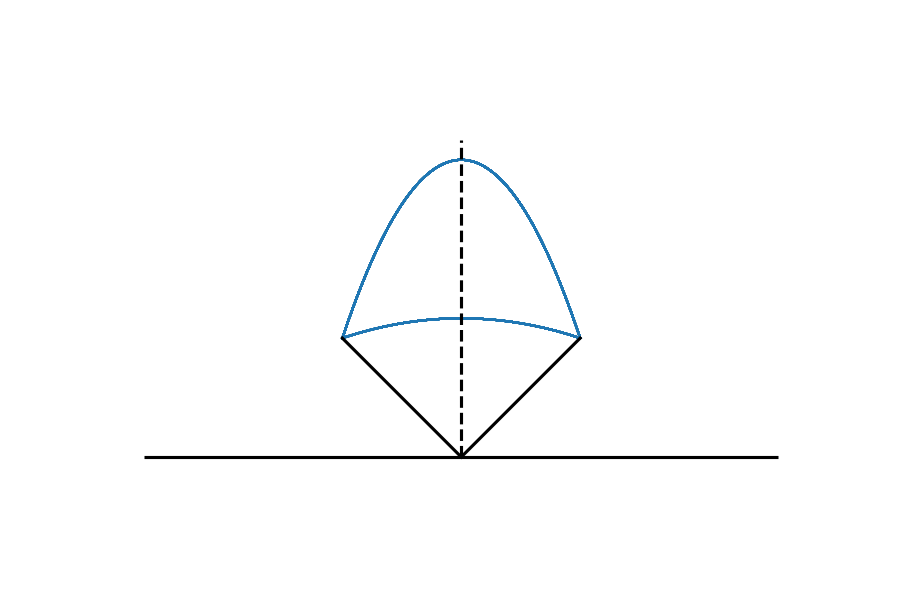}
    }
    \caption{Periodic orbits of the symmetric orthogonal gravitational wedge
      billiard, represented in configuration space $(x,y)$.}
    \label{fig:owb-periodic-orbits}
  \end{center}
\end{figure}

This work generalises the work done by Lehtihet and Miller
\cite{lehtihet1986numerical}, Richter et. al. \cite{richter1990breathing}, and
Szeredi \cite{szeredi1993classical,szeredi1996hard} on the symmetric orthogonal
gravitational wedge billiard.
For that specific model, one can show \cite{anderson2019thesis} the existence of
a period-1 orbit in the billiard, illustrated in Figure
\ref{fig:owb-periodic-orbits}, which corresponds to fixed-point solutions of
the collision maps.

\section{Model}
\label{sec:rowb}

\subsection{Geometry}\label{sec:geometry}
\begin{figure}[ht]
  \begin{center}
    \includegraphics{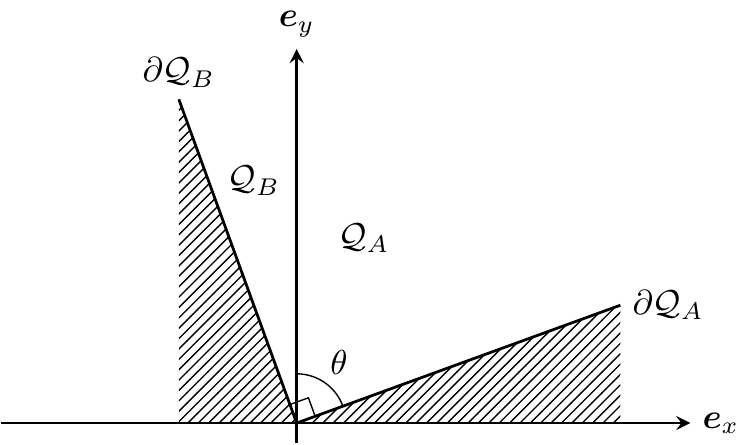}
    \caption{Regions of motion and boundaries for the rotated orthogonal wedge
    billiard.}
    \label{fig:rowb-regions-of-motion}
  \end{center}
\end{figure}

We consider the motion of a point particle of mass $m$ within a constant
gravitational field $\bm{g}$, its motion restricted to the regions
\begin{align*}
  \mathcal{Q}_A &= \left\{ (x,y) \in \mathbb{R}^2 : x \geq 0, \; y \geq
  x\cot(\theta) \right\}, \\
  \mathcal{Q}_B &= \left\{ (x,y) \in \mathbb{R}^2 :
  x < 0, \; y \geq x\tan(\theta) \right\},
\end{align*}
with corresponding boundaries
\begin{align*}
  \partial \mathcal{Q}_A &= \left\{ (x,y) \in \mathbb{R}^2 : x \geq 0, \; y =
  x\cot(\theta) \right\}, \\
  \partial \mathcal{Q}_B &= \left\{ (x,y) \in
  \mathbb{R}^2 : x \leq 0, \; y = x\tan(\theta) \right\},
\end{align*}
as illustrated in Figure \ref{fig:rowb-regions-of-motion}.
The real-valued variable $\theta$ represents the angle measured clockwise from
the vertical and may take values on the interval $(0, \pi/2)$.

We shall call the set $\partial \mathcal{Q} \defeq \partial \mathcal{Q}_A \cup
\partial \mathcal{Q}_B$ the \emph{rotated orthogonal wedge}.
The boundary $\partial \mathcal{Q}_A$ shall be called the \emph{right-hand
slope} or \emph{right-hand wall} of the wedge; similarly, the boundary $\partial
\mathcal{Q}_B$ shall be called the \emph{left-hand slope} or
\emph{left-hand wall} of the wedge.
The intersection of the two boundaries $\partial \mathcal{Q}_A \cap \partial
\mathcal{Q}_B$ shall be called the \emph{wedge vertex}.
The region $\mathcal{Q} \defeq \mathcal{Q}_A \cup \mathcal{Q}_B$ is called the
\emph{region of allowed motion}.

\begin{figure}
  \begin{center}
    \includegraphics{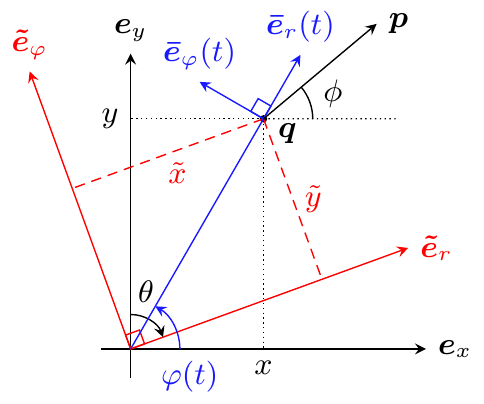}
    \caption{Reference systems used in our study of the rotated orthogonal wedge
    billiard. (Colour online)}
    \label{fig:rowb-radial-ref-sys}
  \end{center}
\end{figure}

We introduce an inertial Cartesian reference system such that the origin $O$ is
fixed at the wedge vertex, the reference axes with unit vectors $\bm{e}_x =
\begin{bmatrix}1 & 0\end{bmatrix}^T$ and $\bm{e}_y = \begin{bmatrix}0 &
1\end{bmatrix}^T$ are orthogonal to each other and directed along the
horizontal and vertical, respectively, as illustrated in Figure
\ref{fig:rowb-radial-ref-sys}.
We shall denote by $\mathcal{B}_C$ the set $\{\bm{e}_x, \bm{e}_y\}$.

Of import in later sections, will be the additional reference systems
$\mathcal{B}_R = \left\{\bm{\bar{e}}_r, \bm{\bar{e}}_\varphi\right\}$ with
origin coinciding with the particle, and $\mathcal{B}_W =
\left\{\bm{\tilde{e}}_r, \bm{\tilde{e}}_{\theta}\right\}$ with origin coinciding
with that of $\mathcal{B}_C$, as illustrated in Figure
\ref{fig:rowb-radial-ref-sys}.
The former, introduced by Lehtihet and Miller \cite{lehtihet1986numerical}, is
dynamic and changes during the motion of the particle; while the
latter, introduced by Szeredi \cite{szeredi1996hard}, is fixed with its
reference vectors coinciding with the wedge walls.


\subsection{Mechanics}
Let $t$ be the variable representing time and let $\bm{q} \defeq \bm{q}(t) \in
\mathcal{Q}$ represent the position vector at some time $t$, and let $\bm{p}
\defeq \bm{p}(t) \in \mathbb{R}^2$ represent the momentum vector of the particle
at some time $t$.
%
In the Cartesian reference system, under a transformation to dimensionless
quantities \cite{anderson2019thesis}, we have the coordinate vectors $[\bm{g}]_C
= \begin{bmatrix}0 & -1\end{bmatrix}^T$, $[\bm{q}]_C = \begin{bmatrix}x &
y\end{bmatrix}^T$ and $[\bm{p}]_C = \begin{bmatrix}u & w\end{bmatrix}^T$.
The corresponding Hamiltonian function is
\begin{equation}
  H(x, y, u, w) = \frac{u^2 + w^2}{2} + y.
  \label{eqn:particle-hamiltonian-dimensionless}
\end{equation}

Integrating the Hamiltonian equations of motion, derived from the Hamiltonian
\eqref{eqn:particle-hamiltonian-dimensionless}, we find that the particle moves
along a parabolic path until it collides with either of the two walls $\partial
\mathcal{Q}_j$ ($j = \{A, B\}$).
These collisions are assumed to be elastic and obey the reflection law.
Since the motion between collisions is completely determined, we focus only on
the collision points themselves.
The collision maps simplify remarkably in the $\mathcal{B}_R$ reference system.
For successive collisions on $\partial \mathcal{Q}_A$, we have the map $F_A :
\partial \mathcal{Q}_A \to \partial \mathcal{Q}_A$ given by
\begin{subequations}\label{eqn:rowb-collision-map-FA}
  \begin{align}
    \bar{u}_{j+1} &= \bar{u}_j - 2\bar{w}_j\cot(\theta),
    \label{eqn:rowb-collision-map-FA-u} \\
    \bar{w}_{j+1}^2 &= \bar{w}_j^2. \label{eqn:rowb-collision-map-FA-w}
  \end{align}
\end{subequations}
Similarly, for successive collisions on $\partial \mathcal{Q}_B$, we have the
map $G_A : \partial \mathcal{Q}_B \to \partial \mathcal{Q}_B$ given by
\begin{subequations}\label{eqn:rowb-collision-map-GA}
  \begin{align}
    \bar{u}_{j+1} &= \bar{u}_j + 2\bar{w}_j\tan(\theta),
    \label{eqn:rowb-collision-map-GA-u} \\
    \bar{w}_{j+1}^2 &= \bar{w}_j^2. \label{eqn:rowb-collision-map-GA-w}
  \end{align}
\end{subequations}
For a collision with the opposite wall $\partial \mathcal{Q}_B$ for the particle
starting on $\partial \mathcal{Q}_A$, we have the map $F_B : \partial
\mathcal{Q}_A \to \partial \mathcal{Q}_B$ given by
\begin{subequations}\label{eqn:rowb-collision-map-FB}
  \begin{align}
    \bar{u}_{j+1} &= \bar{w}_j - \left(\bar{u}_j +
    \bar{w}_{j+1}\right)\tan(\theta), \label{eqn:rowb-collision-map-FB-u} \\
    \bar{w}_{j+1}^2 &= 2E - \bar{w}_j^2. \label{eqn:rowb-collision-map-FB-w}
  \end{align}
\end{subequations}
Similarly, for a collision with $\partial \mathcal{Q}_A$ for the particle
starting on $\partial \mathcal{Q}_B$, we have the map $G_B : \partial
\mathcal{Q}_B \to \partial \mathcal{Q}_A$ given by
\begin{subequations}\label{eqn:rowb-collision-map-GB}
  \begin{align}
    \bar{u}_{j+1} &= -\bar{w}_j - \left(\bar{u}_j -
    \bar{w}_{j+1}\right)\cot(\theta) \label{eqn:rowb-collision-map-GB-u} \\
    \bar{w}_{j+1}^2 &= 2E - \bar{w}_j^2. \label{eqn:rowb-collision-map-GB-w}
  \end{align}
\end{subequations}

\subsection{One-dimensional approximation of motion}
We now consider the transformation of the Hamiltonian from the Cartesian
reference system $\mathcal{B}_C$ to the reference system $\mathcal{B}_W$.
(The transformations between the different reference systems is detailed in
Appendix \ref{sec:coordinate-transforms}.)

Let $\bm{q} = \begin{bmatrix} x & y \end{bmatrix}^T$ and $\bm{p} =
\begin{bmatrix} u & w \end{bmatrix}^T$ denote the coordinate vectors in
$\mathcal{B}_C$, as previously mentioned.
We shall denote the corresponding position and momentum coordinate vectors in
$\mathcal{B}_W$ by $\bm{\tilde{q}} = \begin{bmatrix}\tilde{x} &
\tilde{y}\end{bmatrix}^T$ and $\bm{\tilde{p}} = \begin{bmatrix} \tilde{u} &
\tilde{w} \end{bmatrix}^T$.

Rewriting the Hamilton function \eqref{eqn:particle-hamiltonian-dimensionless}
in terms of the $\mathcal{B}_W$ coordinates yields
\begin{equation}\label{eqn:particle-hamiltonian-dimensionless-radial}
  \tilde{H} = H(\tilde{x}, \tilde{y}, \tilde{u}, \tilde{w}) =
  \frac{\tilde{u}^2}{2} + \tilde{x}\cos(\theta) + \frac{\tilde{w}^2}{2} +
  \tilde{y}\sin(\theta).
\end{equation}
Define
\begin{equation}\label{eqn:hamiltonians-onedim}
  \tilde{H}_{\tilde{x}} \defeq \frac{\tilde{u}^2}{2} + \tilde{x}\cos(\theta),
  \quad
  \tilde{H}_{\tilde{y}} \defeq \frac{\tilde{w}^2}{2} + \tilde{y}\sin(\theta)
\end{equation}
and note that both $H_{\tilde{x}}$ and $H_{\tilde{y}}$ are Hamiltonians
associated with the one dimensional motion of a particle in a (rotated)
gravitational field.
Thus, if we suppose that $\tilde{y}$, $\tilde{w}$ is small enough, which
corresponds to particle motion very close to $\partial \mathcal{Q}_A$, then
$H_{\tilde{y}} \to 0$ and $\tilde{H} \approx H_{\tilde{x}}$.
A similar argument holds for $H_{\tilde{y}}$ for very small $\tilde{x},\
\tilde{u}$.

The equations of motion for the particle very close to $\partial \mathcal{Q}_A$
are
%
\begin{subequations}\label{eqn:particle-hx-eom-solutions}
  \begin{align}
    \tilde{x}(t) &= t\tilde{u}_0 - \frac{t^2\cos(\theta)}{2},
    \label{eqn:particle-hx-eom-solutions-x} \\
    \tilde{u}(t) &= \tilde{u}_0 - t\cos(\varphi)
    \label{eqn:particle-hx-eom-solutions-u}
  \end{align}
\end{subequations}
and for the particle very close to $\partial \mathcal{Q}_B$ are
\begin{subequations}\label{eqn:particle-hy-eom-solutions}
  \begin{align}
    \tilde{y}(t) &= t\tilde{w}_0 - \frac{t^2\sin(\theta)}{2},
    \label{eqn:particle-hy-eom-solutions-y} \\
    \tilde{w}(t) &= \tilde{w}_0 - t\sin(\theta).
    \label{eqn:particle-hy-eom-solutions-w}
  \end{align}
\end{subequations}
From the Hamiltonians \eqref{eqn:hamiltonians-onedim} we may derive bounds on
the trajectories in the configuration space:
\begin{equation}\label{eqn:orwb-trajectory-bounds}
  0 \leq \tilde{x}(t) \leq \frac{E}{\cos(\theta)}, \quad 0 \leq
  \tilde{y}(t) \leq \frac{E}{\sin(\theta)},
\end{equation}
where $E$ is the constant energy fixed at the start of the particle's motion.

\section{Dynamics}
\subsection{Fixed points of the collision maps}
The maps $F_A$ and $G_A$ have the family of fixed points $(\bar{u}^*, \bar{w}^*)
= (c, 0)$ where $c \in \mathbb{R}$ is a constant.
This corresponds to the particle sliding up ($c > 0$) or down ($c < 0$) either
$\partial \mathcal{Q}_A$ or $\partial \mathcal{Q}_B$

For the maps $F_B$ and $G_B$, we obtain
\begin{equation}\label{eqn:fp-fb}
  \bar{u}^* = \sqrt{E}\left(\frac{1 - \tan(\theta)}{1 + \tan(\theta)}\right),
  \qquad
  \bar{w}^* = \sqrt{E},
\end{equation}
and
\begin{equation}\label{eqn:fp-gb}
  \bar{u}^* = \sqrt{E}\left(\frac{\cot(\theta) - 1}{\cot(\theta) + 1}\right),
  \qquad
  \bar{w}^* = \sqrt{E},
\end{equation}
respectively.
In the symmetric wedge billiard, these two fixed points are identical and
correspond to a period one trajectory, as illustrated by Figure
\ref{fig:owb-fixed-point-trajectory}.
To obtain a similar period one trajectory for the rotated orthogonal wedge
billiard, we would need to reflect the momentum component across both
$\bar{\bm{e}}_r$ and $\bar{\bm{e}}_\varphi$ at the collision point, and set it
equal to \eqref{eqn:fp-fb}, which results in
$
\begin{bmatrix}
  \bar{u}' & \bar{w}'
\end{bmatrix}^T
=
\begin{bmatrix}
  -\bar{u} & -\bar{w}
\end{bmatrix}^T.
$
However, such a reflection is only possible when $\theta = \pi/4$.

\subsection{Periodic orbits}
\begin{figure}
  \begin{center}
    \subfloat[$\theta = 50^\circ$]{
      \includegraphics[width=.45\textwidth]{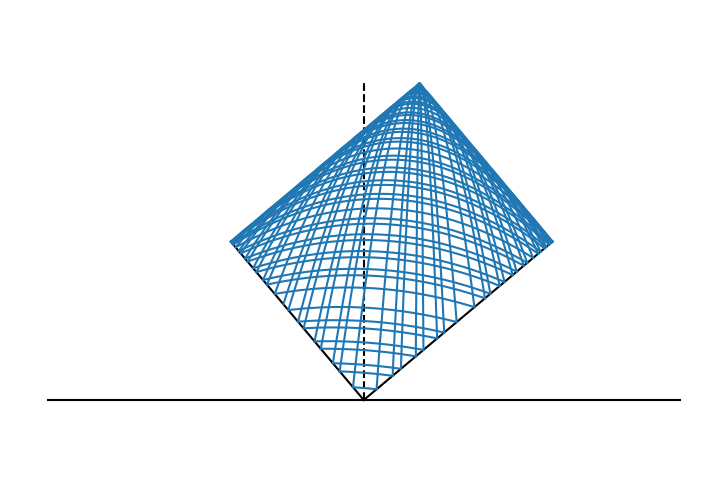}
    }
    \subfloat[$\theta = 60^\circ$]{
      \includegraphics[width=.45\textwidth]{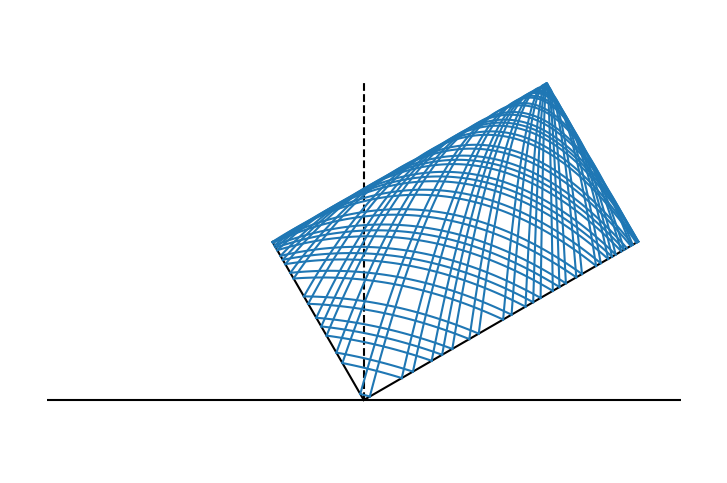}
    }
    \\
    \subfloat[$\theta = 70^\circ$]{
      \includegraphics[width=.45\textwidth]{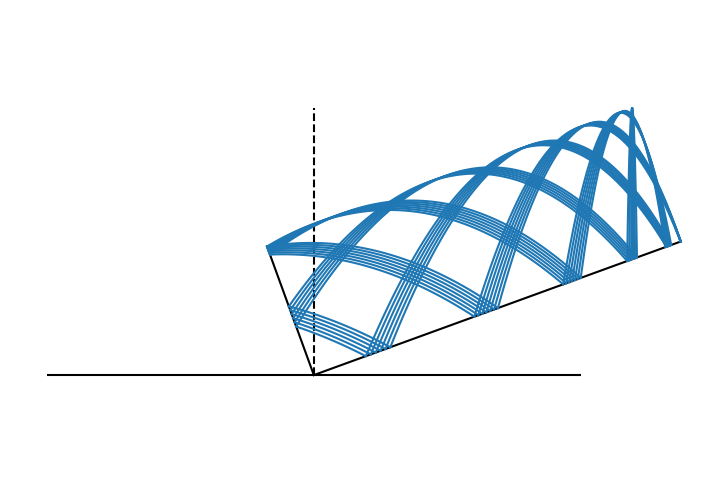}
    }
    \caption{Dense trajectories in the rotated orthogonal wedge billiard for
    various $\theta$, $\phi = 90^\circ$, $\bar{u}_0 = 0$ and $\bar{w}_0 = 1$,
    after $50$ collisions.}
    \label{fig:orwb-trajectory-examples}
  \end{center}
\end{figure}

Using equations \eqref{eqn:particle-hx-eom-solutions-u} and
\eqref{eqn:particle-hy-eom-solutions-w}, we may define the following time maps
for collisions with the wedge walls in the one-dimensional approximation:
\begin{equation}\label{eqn:impact-oscillator-map-x}
  t_{j+1}^A = \frac{2\tilde{u}(t_j)}{\cos(\theta)} =
  2jT_A, \quad
  T_A \defeq \frac{\tilde{u}_0}{\cos(\theta)}
\end{equation}
and
\begin{equation}\label{eqn:impact-oscillator-map-y}
  t_{k+1}^B = \frac{2\bar{w}(t_k)}{\sin(\theta)} = 2kT_B, \quad
  T_B \defeq \frac{\tilde{w}_0}{\sin(\theta)},
\end{equation}
with $j,k \in \left\{0,1,2,\dotsc\right\}$.
The trajectories of the particle's motion will be dense if $T_A/T_B$ is
irrational, illustrated in Figure \ref{fig:orwb-trajectory-examples}.
The ratio $T_A/T_B$ is rational for
\begin{equation}\label{eqn:theta-crit}
  \theta^* = \arctan\left(\frac{p}{q}\right),
\end{equation}
where $p, q \in \mathbb{N}$ such that $p > 0$, $q > 0$ and $\gcd(p,q) = 1$.
To prove this, note that
\[
  \frac{T_A}{T_B} = \frac{\tilde{u}_0}{\tilde{w}_0}\tan(\theta)
\]
which may be simplified to
\begin{equation}\label{eqn:orwb-period-ratio}
  \frac{T_A}{T_B} = \tan(\theta)
\end{equation}
if we make the subsitutions $\tilde{u}_0 = \bar{u}_0$ and $\tilde{w}_0 =
\bar{u}_0$ for the particle close to $\partial \mathcal{Q}_A$ and $\partial
\mathcal{Q}_B$ respectively.
If we substitute $\theta^*$ into equation \eqref{eqn:orwb-period-ratio}, then
\[
  \frac{T_A}{T_B} = \tan(\theta^*) =
  \tan\left(\arctan\left(\frac{p}{q}\right)\right) = \frac{p}{q}
\]
and $T_A/T_B \in \mathbb{Q}$.
The restriction that $p$ and $q$ be positive follows from the restriction on the
allowed values for $\theta$, that is, $\theta \in (0, \pi/2)$.

If $p = q$, then $\theta^* = \pi/4$ which corresponds to the symmetric
orthogonal wedge billiard.
The fixed points of the collision maps \eqref{eqn:rowb-collision-map-FB} and
\eqref{eqn:rowb-collision-map-GB} become identically $(\bar{u}_*, \bar{w}_*^2) =
(0, E)$, which correspond to the period-$1$ orbit in the symmetric orthogonal
wedge billiard, as previously mentioned.

For $p \neq q$, we either have $p > q$, from which follows that $p/q > 1$ and
$\theta^* > \pi/4$, or $p < q$, from which follows that $p/q < 1$ and $\theta^*
< \pi/4$.

\begin{figure}
  \begin{center}
    \subfloat[\label{fig:rowb_periodic_trajectory1} $p = 1$, $q = 2$]{
      \includegraphics[width=.45\textwidth]{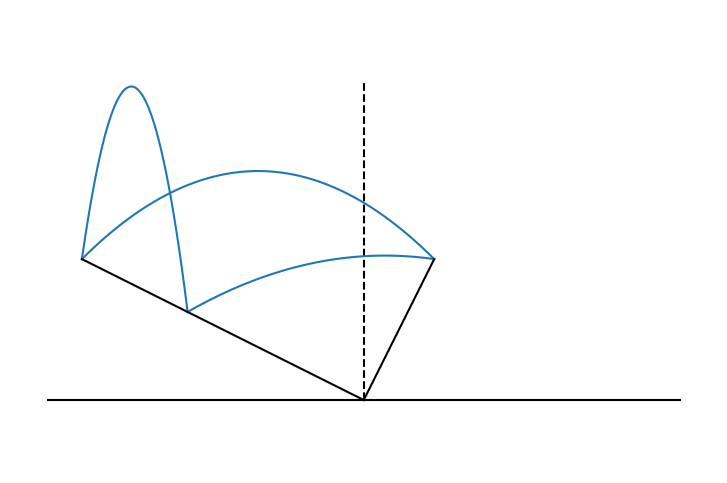}
    }
    \subfloat[\label{fig:rowb_periodic_trajectory3} $p = 1$, $q = 3$]{
      \includegraphics[width=.45\textwidth]{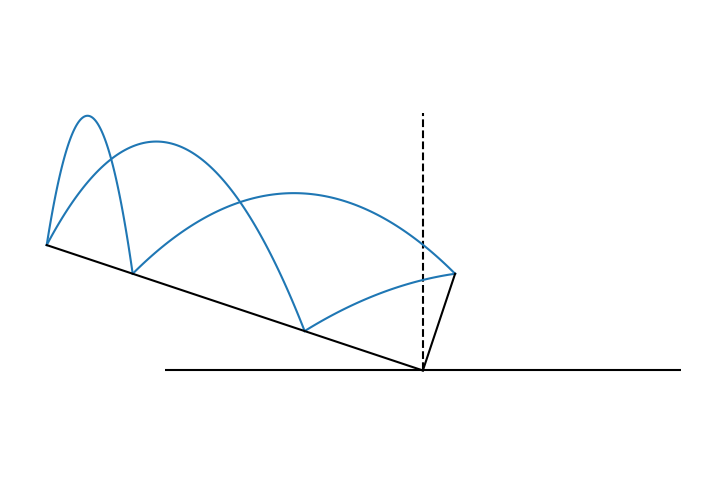}
    }
    \\
    \subfloat[\label{fig:rowb_periodic_trajectory4} $p = 3$, $q = 1$]{
      \includegraphics[scale=0.5]{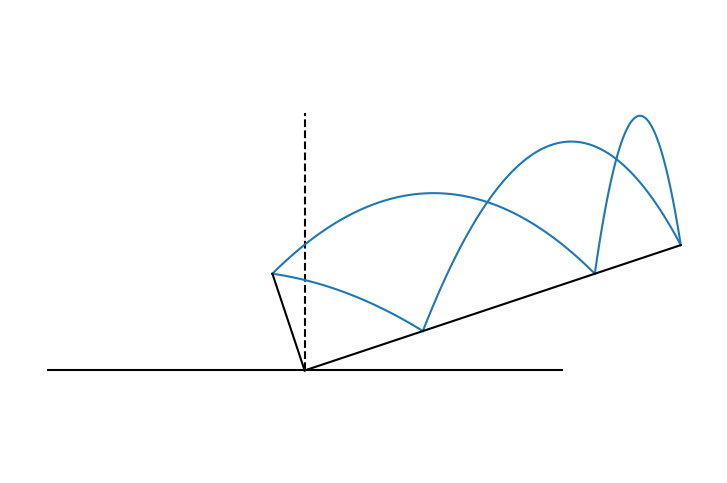}
    }
    \subfloat[\label{fig:rowb_periodic_trajectory5} $p = 2$, $q = 3$]{
      \includegraphics[width=.45\textwidth]{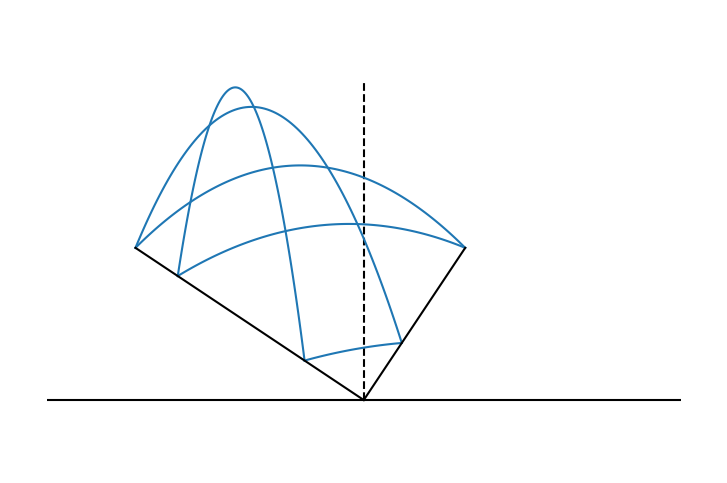}
    }
    \\
    \subfloat[\label{fig:rowb_periodic_trajectory6} $p = 2$, $q = 5$]{
      \includegraphics[width=.45\textwidth]{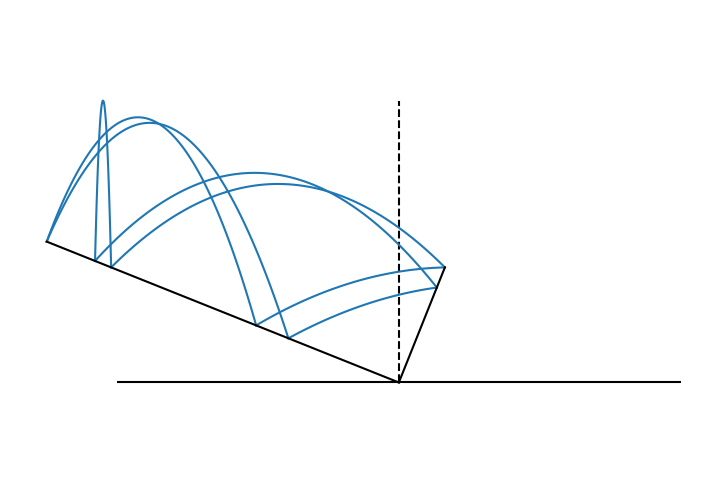}
    }
    \caption{Periodic trajectories for a particle launced from $\partial
    \mathcal{Q}_A$.}
    \label{fig:periodic-orbits}
  \end{center}
\end{figure}

Substituting $\theta^*$ into the fixed point equations \eqref{eqn:fp-fb} and
\eqref{eqn:fp-gb} yields the identical expression
\begin{equation}\label{eqn:rowb-fixed-points-theta-crit}
  \bar{u}^* = \sqrt{E}\left(\frac{q - p}{q + p}\right), \qquad \bar{w}^* =
  \sqrt{E}.
\end{equation}
For a particle starting on $\partial \mathcal{Q}_A$, using
\eqref{eqn:rowb-fixed-points-theta-crit} will yield a periodic orbit of period
$p + q$ with $p$ collisions on $\partial \mathcal{Q}_A$ and $q$ collisions on
$\partial \mathcal{Q}_B$.
The number of collisions per side is determined from
\[
  \frac{T_A}{T_B} = \frac{p}{q} =
  \frac{p\bar{u}_0}{\cos(\theta^*)} \, \frac{\sin(\theta^*)}{q\bar{u}_0} =
  \frac{p\tilde{u}_0}{\cos(\theta^*)} \, \frac{\sin(\theta^*)}{q\tilde{w}_0}
\]
which leads to
\[
  t_{p+1}^A = 2pT_A, \quad t_{q+1}^B = 2qT_B
\]
in equations \eqref{eqn:impact-oscillator-map-x} and
\eqref{eqn:impact-oscillator-map-y}.
Some examples of periodic orbits are illustrated in Figure
\ref{fig:periodic-orbits}.
Numeric simulations show that these periodic orbits are sensitive with respect
to the initial conditions \eqref{eqn:rowb-fixed-points-theta-crit}---a small
perturbation $\bar{u}^* + \varepsilon$ leads to dense orbits once more.

\begin{figure}
  \begin{center}
    \subfloat[\label{fig:periodic-points-phase-space-full}$p,q \in
    \{1,2,3,\dotsc,25\}$]{
      \includegraphics[scale=0.8]{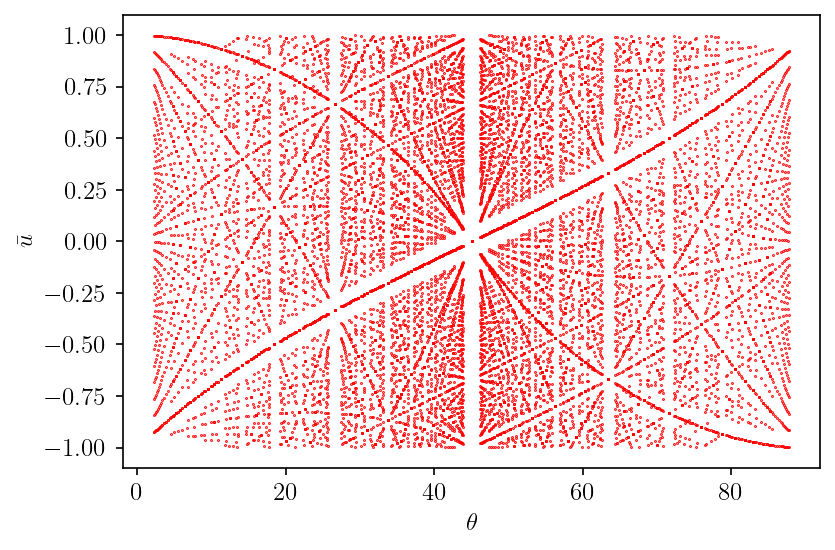}
    }
    \\
    \subfloat[\label{fig:periodic-points-phase-space-half}$p \in
    \{1,2,3,\dotsc,25\},\ q \in \{p + 1, \dotsc, 25\}$]{
      \includegraphics[scale=0.8]{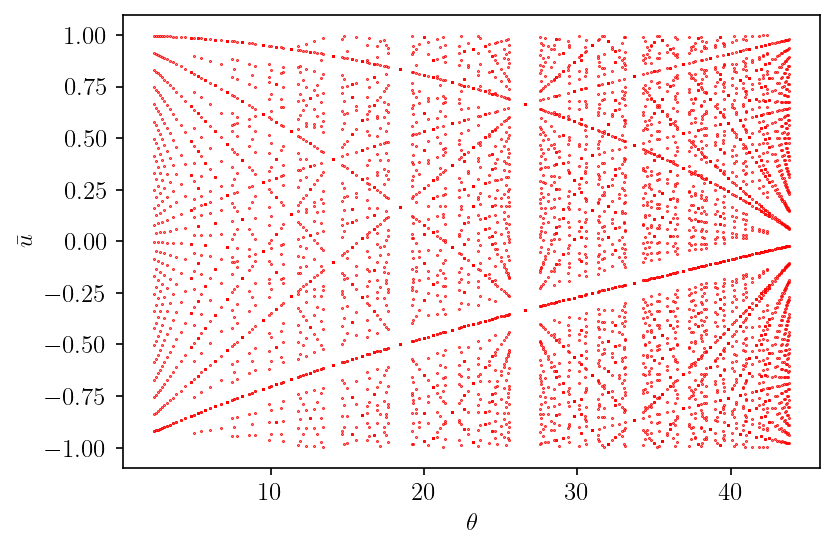}
    }
    \caption{A plot of the $\bar{u}$ component of a periodic trajectory versus
    $\theta$.}
    \label{fig:periodic-points-phase-space}
  \end{center}
\end{figure}

We also plotted the values of $\theta$ versus $\bar{u}$ for various $p$, $q$ in
Figure \ref{fig:periodic-points-phase-space}.
Figure \ref{fig:periodic-points-phase-space-full} indicate a rotational symmetry
in $(\theta, \bar{u})$ space, which allows us to focus on $\theta$-values in
$(0, \pi/4]$ generated from suitable $p$, $q$ values.
Of interest are the ``windows'' which appear for various critical initial
$\theta$, we were unable to determine any special relationship between these
different $\theta$ values with respect each other.

\section{Conclusion}
We studied a special case of the asymmetric gravitational wedge billiard
previously introduced by the author
\cite{anderson2019thesis,anderson2019computational}.
This extends work on the symmetric orthogonal gravitational wedge billiard
previously studied by Lehtihet and Miller \cite{lehtihet1986numerical}, Richter
et. al. \cite{richter1990breathing}, and Szeredi
\cite{szeredi1993classical,szeredi1996hard}.

We derived the conditions under which periodic orbits as well as dense,
non-periodic orbits are found.
For periodic orbits, we derived an explicit formula for the wedge angle and
initial components of the momentum which guarantees a periodic orbit of certain
period and known number of collisions on each wall of the wedge.
However, numeric simulations show that these periodic orbits are unstable with
respect to initial conditions, a small perturbation $\bar{u} + \varepsilon$
leads to dense orbits once more.
We were also unable to reproduce Figure \ref{fig:owb-period-two-trajectory}
using the formula derived for periodic trajectories.
The question is whether this is a degenerate case of the orbit presented in
Figure \ref{fig:owb-fixed-point-trajectory} or not.

\section*{Acknowledgements}
This work extends work done in the author's doctoral thesis
\cite{anderson2019thesis}.
The author wishes to thank the University of Johannesburg for the financial
assistance and opportunity afforded to pursue the doctorate.

\appendix
\section{Coordinate transformations} \label{sec:coordinate-transforms}
In this section, we present the technical details of transformations between the
various reference systems introduced in section \ref{sec:geometry}.
Note the time-dependence of the angle $\varphi$ in Figure
\ref{fig:rowb-radial-ref-sys}---this angle changes as the particle moves,
whereas $\theta$ remains fixed during the motion.
Consider the rotation matrix
\[
  R_\varphi =
  \begin{bmatrix}
    \cos(\varphi) & -\sin(\varphi) \\
    \sin(\varphi) & \cos(\varphi)
  \end{bmatrix}
\]
and note that
\begin{align*}
  \begin{bmatrix}
    \bm{\bar{e}}_r & \bm{\bar{e}}_{\varphi}
  \end{bmatrix}
  &=
  \begin{bmatrix}
    \bm{e}_x & \bm{e}_y
  \end{bmatrix}
  \begin{bmatrix}
    \cos(\varphi) & -\sin(\varphi) \\
    \sin(\varphi) & \cos(\varphi)
  \end{bmatrix}
  =
  \begin{bmatrix}
    \bm{e}_x & \bm{e}_y
  \end{bmatrix}
  R_\varphi, \\
  \begin{bmatrix}
    \bm{\tilde{e}}_r & \bm{\tilde{e}}_{\theta}
  \end{bmatrix}
  &=
  \begin{bmatrix}
    \bm{e}_x & \bm{e}_y
  \end{bmatrix}
  \begin{bmatrix}
    \sin(\theta) & -\cos(\theta) \\
    \cos(\theta) & \sin(\theta)
  \end{bmatrix}
  =
  \begin{bmatrix}
    \bm{e}_x & \bm{e}_y
  \end{bmatrix}
  R_{\pi/2 - \theta}, \\
  \begin{bmatrix}
    \bm{\tilde{e}}_r & \bm{\tilde{e}}_{\theta}
  \end{bmatrix}
  &=
  \begin{bmatrix}
    \bm{\bar{e}}_r & \bm{\bar{e}}_{\varphi}
  \end{bmatrix}
  \begin{bmatrix}
    \sin(\theta + \varphi) & -\cos(\theta + \varphi) \\
    \cos(\theta + \varphi) & \sin(\theta + \varphi)
  \end{bmatrix}
  =
  \begin{bmatrix}
    \bm{\bar{e}}_r & \bm{\bar{e}}_{\varphi}
  \end{bmatrix}
  R_\varphi^TR_{\pi/2 - \theta}
\end{align*}

If we denote the (dimensionless) components of the position and momentum vectors
in the different reference systems analogously, that is,
\begin{gather*}
  [\bm{q}]_C =
  \begin{bmatrix}
    x \\ y
  \end{bmatrix}, \quad
  [\bm{p}]_C =
  \begin{bmatrix}
    u \\ w
  \end{bmatrix}, \quad
  [\bm{q}]_R =
  \begin{bmatrix}
    \bar{x} \\ \bar{y}
  \end{bmatrix}, \quad
  [\bm{p}]_R =
  \begin{bmatrix}
    \bar{u} \\ \bar{w}
  \end{bmatrix}, \\
  [\bm{q}]_W =
  \begin{bmatrix}
    \tilde{x} \\ \tilde{y}
  \end{bmatrix}, \quad
  [\bm{p}]_W =
  \begin{bmatrix}
    \tilde{u} \\ \tilde{y}
  \end{bmatrix},
\end{gather*}
then the transformations of the components between the various reference systems
are
\begin{align*}
  [\bm{q}]_R &= R_\varphi^T [\bm{q}]_C =
  R^T_{\varphi}R_{\pi/2-\theta}[\bm{q}]_W, &
  [\bm{p}]_R &= R_\varphi^T [\bm{p}]_C =
  R^T_{\varphi}R_{\pi/2-\theta}[\bm{p}]_W, \\
  [\bm{q}]_W &= R_{\pi/2 - \theta}^T [\bm{q}]_C =
  \R^T_{\pi/2-\theta}R_{\varphi}[\bm{q}]_R &
  [\bm{p}]_R &= R_{\pi/2 - \theta}^T[\bm{p}]_C =
  R^T_{\pi/2-\theta}R_{\varphi}[\bm{p}]_R.
\end{align*}
Also note that
\[
  [\bm{g}]_R = -
  \begin{bmatrix}
    \sin(\varphi) \\ \cos(\varphi)
  \end{bmatrix},
  \quad
  [\bm{g}]_W = -
  \begin{bmatrix}
    \cos(\theta) \\ \sin(\theta)
  \end{bmatrix}.
\]
For the particle on $\partial \mathcal{Q}_A$, $\varphi = \pi/2 - \theta$ and
\[
  \begin{bmatrix}
    \bm{\tilde{e}}_r & \bm{\tilde{e}}_{\theta}
  \end{bmatrix}
  =
  \begin{bmatrix}
    \bm{\bar{e}}_r & \bm{\bar{e}}_{\varphi}
  \end{bmatrix}
  \begin{bmatrix}
    1 & 0 \\
    0 & 1
  \end{bmatrix}
  =
  \begin{bmatrix}
    \bm{\bar{e}}_r & \bm{\bar{e}}_{\varphi}
  \end{bmatrix}, \qquad
  \begin{bmatrix}
    \tilde{u} \\ \tilde{w}
  \end{bmatrix}
  =
  \begin{bmatrix}
    \bar{u} \\ \bar{w}
  \end{bmatrix}.
\]
For the particle on $\partial \mathcal{Q}_B$, $\varphi = \pi - \theta$ and
\[
  \begin{bmatrix}
    \bm{\tilde{e}}_r & \bm{\tilde{e}}_{\theta}
  \end{bmatrix}
  =
  \begin{bmatrix}
    \bm{\bar{e}}_r & \bm{\bar{e}}_{\varphi}
  \end{bmatrix}
  \begin{bmatrix}
    0 & 1 \\
    -1 & 0
  \end{bmatrix}
  =
  \begin{bmatrix}
    -\bm{\bar{e}}_{\varphi} & \bm{\bar{e}}_r
  \end{bmatrix}, \qquad
  \begin{bmatrix}
    \tilde{u} \\ \tilde{w}
  \end{bmatrix}
  =
  \begin{bmatrix}
    -\tilde{w} \\ \tilde{u}
  \end{bmatrix}.
\]

\bibliographystyle{abbrv}
\bibliography{manuscript_rowb_arxiv}{}
\end{document}